\documentclass{proc-l}

\usepackage{amsthm}
\usepackage{amsmath}
\usepackage{amsfonts}

\title{Complete Intersections in Toric Ideals}

\author[E. Cattani, R. Curran, and A. Dickenstein]
{Eduardo Cattani, Raymond Curran, and Alicia Dickenstein}
\address{Eduardo Cattani: Department of Mathematics
and Statistics. University
of Massachusetts. Amherst, MA 01003, USA}
\email{cattani@math.umass.edu}
\thanks{E. Cattani is partially supported by NSF Grant
DMS--0099707}

\address{Raymond Curran: Department of Mathematics
and Statistics. University
of Massachusetts. Amherst, MA 01003, USA}
\curraddr{Department of Mathematical and Computer Sciences.
Metropolitan State College of Denver.
Denver, CO  80202, USA. }
\email{rcurran@mscd.edu}
\address{Alicia Dickenstein: Departamento~de
Matematica, FCEyN.
Universidad de Buenos Aires. (1428) Buenos Aires,
Argentina}
\email{alidick@dm.uba.ar}
\thanks{A. Dickenstein is partially supported by
UBACYT X042,
Argentina.}

\def\cocoa{{\hbox{\rm C\kern-.13em o\kern-.07em C\kern-.13em 
o\kern-.15em A}}}

\newcommand{\baseRing}[1]{\ensuremath{\mathbb{#1}}}
\newcommand{\Z}{\baseRing{Z}}
\newcommand{\R}{\baseRing{R}}
\newcommand{\C}{\baseRing{C}}
\newcommand{\N}{\baseRing{N}}
\newcommand{\Q}{\baseRing{Q}}

\theoremstyle{plain}
\newtheorem{theorem}{Theorem}[section]
\newtheorem{lemma}[theorem]{Lemma}
\newtheorem{corollary}[theorem]{Corollary}

\newtheorem{proposition}[theorem]{Proposition}

\theoremstyle{definition}

\newtheorem{definition}[theorem]{Definition}
\newtheorem{remark}[theorem]{Remark}

\newtheorem{example}[theorem]{Example}

\numberwithin{equation}{section}

\newcommand{\Script}[1]{\ensuremath{{\mathcal{#1}}}}

\newcommand{\LL}{\Script{L}}

\begin{document}

\begin{abstract}
We present examples which show that in dimension higher than one
or codimension higher than two, there exist toric ideals $I_A$ such that 
no binomial ideal contained in $I_A$ and of the same dimension is a complete
intersection.  This result has important implications in sparse
elimination theory and in the study of the Horn system of partial
differential equations.
\end{abstract}

\footnotetext[1]{AMS Subject Classification:
Primary 14M10, Secondary
14M25, 13C40}

\maketitle

\section{Introduction}
Given
a configuration
$A = \{a_1,\dots,a_n\}\subset \Z^m$ of integral
points generating $\Z^m$, the (toric) ideal $I_A \subset 
\C[x_1,\dots,x_n]$ is
generated by all binomials
$$x^u - x^v\, ,$$
whose exponents $u,v\in \N^n$ satisfy $A\cdot u = A\cdot v$.
Here we are also denoting by
$A$ the $m\times n$ matrix whose
$j$-th column is $a_j$.
Note that $I_A$ is weighted-homogeneous
for every weight $w$ in the row-span of $A$. We shall assume
that $(1,\dots,1)$ is in the row-span of $A$ and, consequently, that
$I_A$ is a (standard) homogeneous ideal. 

The configuration $A$ and its associated toric ideal $I_A$ are
the central characters in several areas of very active
research in
commutative algebra, algebraic geometry, and differential
algebra. The associated projective variety $X_A$ has a
natural action of the algebraic torus $(\C^*)^m$ making
it into a projective toric variety of dimension
$d:=m-1$. We call $d$ the dimension of $A$.

Let $\LL_A \subset \R^n$ be the
lattice
$$ \LL_A = \{v \in \Z^n : A\cdot v = 0\}$$
and let $r:=n-m$ be the codimension of $A$.
Let  $B= \{v_1,\dots,v_r\} \subset \LL_A$ 
be a maximal set of vectors linearly independent
over $\Q$.
We also denote by
$B$ the $n\times r$ matrix whose
$j$-th column is $v_j$ and
consider the
binomial ideal $J_B$ generated by the binomials
$x^{v_j^+} - x^{v_j^-}$, where $v_j = {v_j^+} - {v_j^+}$ is
the decomposition in positive and negative components.

It has been proven in \cite{es} that $J_B \otimes \C[x_1^{\pm1},
\dots, x_n^{\pm 1}]$ is  always a complete intersection
in the Laurent polynomial ring, and that
it coincides with $I_A \otimes \C[x_1^{\pm1},
\dots, x_n^{\pm 1}]$ if and only if the greatest common divisor $g$
of the maximal minors of $B$ satisfies $g=1$. This positive
integer $g$ is precisely the index with respect to  $\LL_A$
of the lattice 
spanned by $B$. 
We will refer to $J_B$ as a {\it  basis ideal\/}.
Following \cite{hs} we will reserve the term
 lattice basis ideal for the case when $B$ is
 a $\Z$-basis of $\LL_A$.
Clearly $J_B \subset I_A$ and, in general, this containment
is proper.

In this note we study 
a question which arises naturally 
in the study of $A$-discriminants \cite{ds}
and of Horn systems of differential equations \cite{dms}:
does every 
 toric ideal contain a
complete intersection  basis ideal (i.e. a complete intersection
binomial ideal of the same dimension)? This is indeed
the case if $X_A$ is a monomial curve ($d=1$) or in the codimension 
two case ($r=2$). The purpose of this note
is to show that in any dimension higher than one and any
 codimension higher than two, there exist toric ideals $I_A$ such that 
no  basis ideal contained in $I_A$ is a complete
intersection. 

We recall that the dual variety $X_A^*$ of $X_A$ is
the Zarisky closure of the locus of hyperplanes tangent to $X_A$
at a smooth point.
When $X_A^*$ is
a hypersurface, its defining equation is the $A$-discriminant.
This notion, introduced by
Gel'fand, Kapranov, and Zelevinsky \cite{gkzbook}, 
generalizes the
classical notion of the discriminant of a univariate polynomial.
Dickenstein and Sturmfels \cite{ds} have shown how
to compute $A$-discriminants
for codimension two (i.e. $n - m =2$) configurations.
A key ingredient of this work is the fact that every 
basis ideal in a codimension two toric ideal is a complete
intersection. Our results show that a different approach
is needed  to describe 
$A$-discriminants in higher codimensions.

The work of Gel'fand et al. on sparse elimination was a step
toward the study of $A$-hypergeometric (or GKZ) systems. Consider
the ``quantized" version of the ideal $I_A$, that is, the left
ideal $H_A(\beta)$ in the Weyl algebra
$$D_n = \C\langle x_1,\dots,x_n,\partial_1,\dots,\partial_n\rangle$$
generated by the toric operators $\partial^u - \partial^v$,
$A\cdot u = A\cdot v$ together with the Euler operators
associated with the $(\C^*)^m$ action:
$$ a_{j1} x_1 \partial_1 + \cdots + a_{jn} x_n \partial_n - \beta_j\,,$$
where $\beta = (\beta_1,\dots,\beta_m)\in \C^d$. An
$A$-hypergeometric function of degree $\beta$ is a locally defined
(multivalued) holomorphic function $\C^n$ annihilated by $H_A(\beta)$.
This notion of hypergeometric functions encompasses most
of the classical univariate and multivariate hypergeometric
functions. GKZ systems are holonomic for all choices of parameters
$\beta$, and so in particular, the corresponding spaces of $A$-hypergeometric
functions are finite dimensional. 
Another classical multivariable generalization of hypergeometric
differential equations is given by  the Horn systems, which are
closely related to the GKZ systems.  A Horn system consists of 
the Euler operators  and
only those toric operators coming from a  basis ideal $J_B$
contained in $I_A$. Dickenstein, Matusevich,
and Sadykov \cite{dms} have shown that, in codimension two, the behavior of an
$A$-hypergeometric
and that of any of its associated Horn systems  is not very different.  
This relies on the fact that
such  basis ideals are complete intersections.
However, if the basis ideal is not a complete
intersection, the Horn system has an infinite dimensional local
solution space for all choices of parameters
and, thus, is never holonomic. Our examples show that
in the general multivariate case, the $A$-hypergeometric and Horn systems
have essentially distinct behavior.

\smallskip

\noindent{\bf Acknowledgments:} We are grateful to
Bernd Sturmfels for many helpful conversations and to
John Abbott and Lorenzo Robbiano of the \cocoa\ group for their programming assistance.  We would also like to thank an anonymous referee for
very useful comments.

\section{Preliminaries}
The study of binomial ideals is intimately connected with the
study of (affine) semigroup algebras. It is from this perspective,
and beginning with the work 
of Herzog \cite{herzog} and Delorme \cite{delorme}, that
the question of classifying complete intersection
binomial ideals has been extensively studied by many authors
\cite{bmt1, bmt2, bmt3, fms1, fms2, fs,hs, nakajima, rosales, sss, 
stanley,
thoma1, thoma2, watanabe}. A combinatorial characterization of these
ideals is given in \cite{fs} in terms of a choice of $B$  
and the notion of mixed matrices; that is, matrices such that every
column contains a strictly positive and a strictly negative
entry. Note that since the columns of the matrix 
$B$ add up to zero, $B$ is automatically mixed.
The following
result follows from \cite[Theorem~2.1]{es} and \cite[Theorem~2.3]{fs}
(see also \cite[Theorem~2.7]{sss}):

\begin{theorem}\label{theo:fs}
The ideal $J_B$ is a complete intersection if and only if
for every mixed $n'\times r'$-submatrix $B'\subset B$ we
have $n'\geq r'$.
\end{theorem}

Since a mixed submatrix must contain at least two rows
it follows that:

\begin{corollary}
If $r\leq 2$, every  basis ideal $J_B$ is a complete
intersection.
\end{corollary}

It is also easy to prove that if $m=2$, i.e. for $X_A$ 
a monomial projective curve, there exists $B$ such
that $J_B$ is a complete intersection. Indeed, let
\begin{equation}\label{curve}
A\ =\
\left(
\begin{array}{cccc}
1 & 1 & \cdots & 1 \\
a_1 & a_2 &\cdots & a_n
\end{array}
\right)\,,
\end{equation}
where $a_1\leq\cdots\leq a_n$ are coprime. 
Performing a row operation that does not change $I_A$, we may
assume without loss of generality that $a_1=0$ and,
therefore all $a_j\geq 0$.
Consider now the following choice $B$ of $A$:
\begin{equation}\label{galecurve}
B\ =\
\left(
\begin{array}{ccccc}
a_3-a_2 & a_4-a_3 & \cdots & a_n-a_{n-1} \\
-a_3 & 0 &\cdots & 0\\
a_2 & -a_4 &\cdots & 0\\
0 & a_3 &\ddots & 0\\
\vdots & \vdots & \ddots &-a_n\\
0 & 0 &\cdots & a_{n-1}\\
\end{array}
\right)\,,
\end{equation}
Since $a_i-a_{i+1}\geq 0$, it follows that every mixed submatrix
of $B$ must contain more rows than columns and therefore
$J_B$ is a complete intersection.

Note that  already in the simplest case of
the twisted cubic; i.e. the curve $X_A$ associated with:
\begin{equation*}
A\ =\
\left(
\begin{array}{cccc}
1 & 1 & 1 & 1 \\
0 & 1 &2 & 3
\end{array}
\right)\,,
\end{equation*}
the toric ideal $I_A$ is not a complete intesection but it does
contain a complete intersection basis ideal.

\bigskip

Given a vector $v\in \Z^n$ we define the support of $v$:
$${\rm supp}(v) \ :=\ \{i\in\{1,\dots,n\}: v_i \not= 0\}.$$
Similarly we set ${\rm supp}^+(v) \ :=\ \{i: v_i > 0\}$,
$\ {\rm supp}^-(v) \ :=\ \{i: v_i < 0\}$. If
$u,v\in \Z^n$ we say that $u$ is {\em conformal} to $v$ if
${\rm supp}^+(u)\subset {\rm supp}^+(v)$ and
${\rm supp}^-(u)\subset {\rm supp}^-(v)$.

Given a configuration $A = \{a_1,\dots,a_n\}\subset\Z^m$, a vector 
$v\in \LL_A$ is
called a {\em circuit} if its support is minimal among
all elements in $\LL_A$ relative to inclusion.
As shown in the proof of \cite[Lemma~4.9]{sturmfels},
if $u$ is
a circuit in $\LL_A$ of maximal dimension $m$ with 
${\rm supp}(u) = \{i_1, \dots, i_{m+1} \}$, then up to multiple:
\begin{equation}\label{circuit}
u\ =\ \sum_{j=1}^{m+1}\,(-1)^j\,\det\left(a_{i_i},\dots,a_{i_{j-1}},
a_{i_{j+1}},\dots,a_{i_{m+1}}\right)\,e_{i_j}.
\end{equation}
The following is Lemma~4.10 in \cite{sturmfels}:

\begin{lemma}\label{lemma:circuit}
Every vector $v \in \LL_A$ may be written as a
non-negative rational linear
combination of $n-m$ circuits each of which is conformal
to $v$.
\end{lemma}

\begin{definition}
Given a configuration $A = \{a_1,\dots,a_n\}\subset\Z^m$
we say that a  basis ideal $J_B$ is generated by
circuits if and only if each of the column vectors of $B$
is a circuit in $\LL_A$.
\end{definition}

\begin{proposition}\label{prop:circuits}
A toric ideal $I_A$ contains a complete
intersection  basis ideal if and only if it contains
a complete intersection  basis ideal generated by circuits.
\end{proposition}

\begin{proof}
Let $\{v_1,\dots,v_r\}\subset\LL_A$
be a $\Q$-linearly independent set
defining a complete intersection  basis
ideal. Let $i\in \{1,\dots,r\}$ be the smallest
index such that
$v_i$ is not a circuit. Using Lemma~\ref{lemma:circuit}
write $v_i$ as a non-negative rational combination of
circuits $w^i_1,\dots,w^i_r$:
$$v_i = q_1 w^i_1 + \cdots + q_r w^i_r\,.$$
Clearly, some $w^i_j$ must be linearly independent
from $\{v_k, k\not=i\}$, hence we may replace
$v_i$ by the circuit $w^i_j$ to obtain a new $\Q$-linearly
independent set. Continuing in this manner we obtain
a linearly independent set $\{w_1,\dots,w_r\}$ consisting of circuits
and such that $w_i$ is conformal to $v_i$.

Let $B$ (respectively $C$) denote the matrix whose
columns are the vectors $v_i$ (respectively $w_i$).
Since $w_i$ is conformal to $v_i$ it follows that
if $C'$ is a mixed submatrix of $C$ then the corresponding
submatrix $B'$ of $B$ is also mixed. Hence, by Theorem~\ref{theo:fs},
if $J_B$ is a complete intersection, so is $J_C$.
\end{proof}

\section{Toric ideals of dimension at least two}

In this section we will exhibit examples of  configurations in
any dimension $d = m-1 \geq 2$,
which do not contain any complete intersection  basis
ideal.
By assumption we may suppose that $A \subset \{1\}\times \Z^{d}$.
By abuse of notation we will
identify $A$ with its projection onto $\Z^{d}$. 

\begin{theorem}\label{theo:dim2}
Let $A$ be
the vertex set of a lattice polytope in  $\R^{d}$
and  $n = |A|$.  There
exists $N = N(d)$ such that, for $n\geq N$, $I_A$ does not contain
any complete intersection  basis ideal.
\end{theorem}

\begin{proof} Let $B \subset \Z^{n \times (n-m)}$  such that  its columns 
 are a $\Q$-basis of the kernel lattice $\LL_A$. Note that every column of
$B$ must contain at least two strictly positive and two strictly 
negative
entries since, otherwise, one of the points in $A$ would be in the
convex hull of some of the other points in $A$ and that is impossible
by assumption.

Suppose $J_B$ is a complete intersection. It then follows
from Theorem~\ref{theo:fs} that every $(n-m-1)\times(n-m)$ minor
of $B$ must be non-mixed. In other words, for any subset
$J\subset \{1,\dots,n\}$, $|J|=m+1$, there exists a column
$v_i$ of $B$ such that either the positive support
${\rm supp}^+(v_i)$
or the negative support ${\rm supp}^-(v_i)$ is contained in $J$.
Clearly, there are
${n \choose d+2}$ index sets $J$ of cardinality $m+1 = d+2$.

Now, for a given $v_i$, its positive support contains at
least two indices and therefore it may be contained in at most
${{n-d}\choose{ {2}}}$ distinct four-index sets $J$.
Thus, since $B$ has $n-d-1$ columns and taking into account the
positive and negative supports, the condition may be satisfied
for at most
$$ 2 \, (n-d-1) \,{n-2 \choose d}$$
index sets $J$. But for $n$ sufficiently
large
\begin{equation}\label{ineq}
{n \choose d+2}\ > 2 \, (n-d-1) \,{n-2 \choose d},
\end{equation}
since the left hand side is a polyomial in $n$ of degree $d+2$
with positive leading term, while the right hand side is a polynomial
in $n$ of degree $d+1$.
\end{proof}

\begin{remark}\label{remark:plane} In the planar case $d=2$,
 the inequality (\ref{ineq}) is satisfied for $n\geq 22$. 
However, it is clear that the estimates
above are very rough and that one should expect Theorem~\ref{theo:dim2}
to hold for $n$ considerably smaller than $22$.

\end{remark}

\begin{example}
Consider a configuration $\{a_1,\dots,a_{10}\}$ of ten points in $\Z^2$
which are the vertices of a polygon.  We may assume
them to be ordered
counterclockwise. Given four indices
$1\leq i<j<k<\ell\leq 10$, there exists a relation
$$\lambda_i a_i - \lambda_j a_j + \lambda_k a_k - \lambda_\ell a_\ell 
=0,$$
where $\lambda_i, \lambda_j, \lambda_k, \lambda_\ell$ are positive
integers. Such a relation defines a circuit in the lattice
kernel $\LL_A$. Using the computer algebra system \cocoa\  \cite{cocoa},
we searched for sets of seven such relations satisfying the condition
in Theorem~\ref{theo:fs}. The following is such an example. We have
only indicated the sign of the coefficients since that is all that
matters in Theorem~\ref{theo:fs} and, for generic coefficients, the
matrix $B$ will be of maximal rank. 

$$B = \left(
\begin{array}{ccccccc}
+ & + & + & 0 & + & 0 & 0 \\
- & - & 0 & 0 & 0 & 0 & 0 \\
+ & 0 & 0 & 0 & - & + & 0 \\
- & + & - & 0 & 0 & 0 & 0 \\
0 & 0 & 0 & + & 0 & - & 0 \\
0 & 0 & 0 & - & + & + & + \\
0 & 0 & 0 & + & 0 & 0 & - \\
0 & 0 & + & 0 & - & 0 & + \\
0 & - & - & 0 & 0 & 0 & 0 \\
0 & 0 & 0 & - & 0 & - & -
\end{array}
\right)
$$
We were unable to obtain similar examples with $n=11$. We
suspect that Theorem~\ref{theo:dim2} holds for polygons in
the plane with at least $11$ vertices.
\end{example}

\section{Toric ideals of codimension at least three}

It follows from Remark~\ref{remark:plane} that 
Theorem~\ref{theo:dim2} furnishes examples of toric
ideals which do not contain any complete intersection
basis ideal in codimension greater than twenty-one.
In this section we will describe a different class of examples
which show that in every
codimension greater than two 
there exist toric ideals
with that same property.

Let $n = m + r$, and consider 
a configuration consisting of the vertices of
a cyclic polytope (we refer to 
\cite{ziegler} for other properties
of this important class of polytopes):
\begin{equation}\label{normal}
A\ =\ \left(
\begin{array}{llll}
1 & 1 & \cdots & 1 \\
t_1 & t_2 & \cdots & t_n \\
t_1^2 & t_2^2 & \cdots & t_n^2 \\
\vdots & \vdots &&\vdots\\
t_1^{m-1} & t_2^{m-1} & \cdots & t_n^{m-1}
\end{array}
\right),
\end{equation}
where 
$0<t_1< t_2 < \cdots<t_n$ are integers. For appropriate choices of $t_1,\dots,t_n$, the columns of $A$ 
span $\Z^m$.

\begin{theorem}\label{theo:codim3}
Let $r\geq 3$.
For $\,n \geq 2(r^2 -r +1) $  the toric ideal $I_A$ associated
with the matrix (\ref{normal}) does not contain any
complete intersection  basis ideal.
\end{theorem}

\begin{proof}

According to Proposition~\ref{prop:circuits}, it suffices to
show that $I_A$ does not contain any complete intersection
 basis ideal generated by circuits which, up to constant, are
 given by the expression (\ref{circuit}).  On the other hand,
note that all maximal
minors of $A$ are non-zero and
the determinant
$$\det\left(a_{i_i},\dots,a_{i_{j-1}},
a_{i_{j+1}},\dots,a_{i_{m+1}}\right)$$
is strictly positive since it
is the Vandermonde determinant for 
$$t_{i_i} < \dots < t_{i_{j-1}} < t_{i_{j+1}} <\dots <t_{i_{m+1}}.$$  
Therefore, if $v\in \LL_A$
is a circuit then it will have exactly $r-1$ zero entries
while the remaining entries will alternate in sign.

This means that if $B$ is an $n\times r$ matrix whose columns are
a circuit basis of $\LL_A$ then each column of $B$ contains
exactly $r-1$ entries which are zero and
$B$ contains a total of $ r(r-1)$ zero entries.
This implies that if $n\geq 2(r^2 -r +1)$, $B$ will
have to contain two consecutive rows all of whose entries
are non-zero. Let $B'$ denote the $2\times 3$ submatrix
consisting of those two rows. Since the signs of the columns
are alternating, $B'$ is mixed and, by Theorem~\ref{theo:fs},
$J_B$ is not a complete intersection. In fact, if 
$B'$ consists of the $j$-th and $k$-th rows of $B$, 
the codimension two ideal $\langle x_j, x_k \rangle$
 is an associated prime ideal of $J_B$.
\end{proof}

\begin{remark}
For $r=3$ the lower bound in Theorem~\ref{theo:codim3} is $14$. One can show that for $n\leq 11$
every toric ideal $I_A$ admits a complete intersection
 basis ideal. We do not know if there are examples
of configurations with $n=12$ or $n=13$ that do not admit any
such  basis ideals.
\end{remark}

\end{document}